\documentclass[a4paper,10pt]{article}
\usepackage{amssymb,amsmath}

\newcommand{\ep}{\varepsilon}

\begin{document}
\title{Corrigendum \\
Quadratic Class Numbers Divisible by 3\\
Functiones et Approximatio, 37 . pp. 203-211.}
\author{D.R. Heath-Brown\\Mathematical Institute, Oxford}
\date{}
\maketitle

In attempting to handle $N_+(X)$ the paper states that the Scholz
reflection principle ``yields $3|h(k)$ for any positive integer for
which $3|h(-3k)$''.  This is not correct, and one cannot establish a
result for $N_+(X)$ in this way.  However one may use a criterion of
Honda \cite[Proposition 10]{H}, which shows that if 
\[27n^2+du^2=4m^3\]
with positive integers $n,u,m,d$, then $3\mid h(d)$ providing that
$(m,3n)=1$ and the polynomial $X^3-mX+n$ has no integer root. This
latter condition is always satisfied if $3\mid m-1$ and $3\nmid n$,
for example.  An
argument completely analogous to that used in the paper then recovers
the stated bound $N_+(X)\gg_{\ep}X^{9/10-\ep}$.

\bigskip
\bigskip

Mathematical Institute,

24--29, St. Giles',

Oxford

OX1 3LB

UK
\bigskip

{\tt rhb@maths.ox.ac.uk}

\end{document}